\documentclass[12pt]{elsarticle}
\usepackage{amsmath, amsthm,amssymb,
MnSymbol,extsizes}

\usepackage[all]{xy}
\usepackage[active]{srcltx}
\usepackage{mathdots}

\newtheorem{theorem}{Theorem}

\newtheorem{lemma}{Lemma}

\newcommand{\od}{^{\oslash}}

\newcommand{\TT}{\mathbb T}
\newcommand{\HH}{\mathbb H}
\newcommand{\CC}{\mathbb C}
\newcommand{\RR}{\mathbb R}
\newcommand{\FF}{\mathbb F}
\newcommand{\mc}{\mathcal}
\DeclareMathOperator{\ind}{ind}
\DeclareMathOperator{\End}{End}

\newcommand{\as}%
{\text{\raisebox{0.75pt}
{$\scriptstyle\medstar$\!}}}
\newcommand{\T}{\top}

\newcommand{\ci}{
\begin{picture}(6,6)
\put(3,3){\circle*{3}}
\end{picture}}

\newcommand{\matt}[1]{\left[\begin{smallmatrix}
   #1\end{smallmatrix}\right]}
\newcommand{\mat}[1]{\begin{bmatrix}
   #1\end{bmatrix}}

\newcommand{\mattc}[1]{
\left(\begin{smallmatrix}#1
\end{smallmatrix}\right)}

\newcommand{\zzz}
{\matt{0&&1\\[-1mm]&\iddots\\1&&0}}

\newcommand{\sykl}[9]
    {\xymatrix@R=25pt@C=45pt{
{\scriptstyle #1}\ar@{->}[d]_{#3}
 \save !<#6 mm,0cm>
\ar@{-}@(ur,dr)^{#4}_{#8\,}
\restore
                      \\
{\scriptstyle #2}
 \save !<#7 pt,0cm>
\ar@{-}@(ur,dr)^{#5}_{#9\,}
\restore }}

\newcommand{\syk}[9]
    {\xymatrix@R=20pt@C=45pt{
{\scriptstyle #1}\ar@{->}[d]_{#3}
 \save !<#6 mm,0cm>
\ar@{-}@(ur,dr)^{#4}_{#8\,}
\restore
                      \\
{\scriptstyle #2}
 \save !<#7 pt,0cm>
\ar@{-}@(ur,dr)^{#5}_{#9\,}
\restore }}

\newcommand{\syrk}[8]
    {\xymatrix@R=20pt@C=45pt{
{\scriptstyle #1}
\ar@{->}[d]_{#3}^{#8}
 \save !<#6 mm,0cm>
\ar@{-}@(ur,dr)^{#4}_{\varepsilon\,}
\restore
                      \\
{\scriptstyle #2}
 \save !<#7 pt,0cm>
\ar@{-}@(ur,dr)^{#5}_{\delta\,}
\restore }}

\newcommand{\sy}[7]
    {\xymatrix@R=20pt@C=45pt{
{\scriptstyle #1}\ar@{->}[d]_{#3}
 \save !<#6 mm,0cm>
\ar@{-}@(ur,dr)^{#4}
\restore
                      \\
{\scriptstyle #2}
 \save !<#7 pt,0cm>
\ar@{-}@(ur,dr)^{#5}
\restore }}

\newcommand{\syyk}[9]
    {\xymatrix@R=30pt@C=45pt{
{\scriptstyle #1}\ar@{->}[d]_{#3}
 \save !<#6 pt,0cm>
\ar@{-}@(ur,dr)^{#4}_{#8\,}
\restore
                      \\
{\scriptstyle #2}
 \save !<#7 pt,0cm>
\ar@{-}@(ur,dr)^{#5}_{#9\,}
\restore }}

\newcommand{\syma}[8]
    {\xymatrix@C=25pt@R=15pt{
{\scriptstyle #1}\ar[d]_{#5}\ar[r]^{ #7}
&{\scriptstyle #2}
                      \\
{\scriptstyle #3}\ar[r]^{#8}&{\scriptstyle #4}\ar[u]_{#6}}}

\newcommand{\symaa}[8]
    {\xymatrix@C=35pt@R=15pt{
{\scriptstyle #1}\ar[d]_{#5}\ar[r]^{ #7}
&{\scriptstyle #2}
                      \\
{\scriptstyle #3}\ar[r]^{#8}&{\scriptstyle #4}\ar[u]_{#6}}}

\newcommand{\symaaa}[8]
    {\xymatrix@C=50pt@R=25pt{
{\scriptstyle #1}\ar[d]_{#5}\ar[r]^{ #7}
&{\scriptstyle #2}
                      \\
{\scriptstyle #3}\ar[r]^{#8}&{\scriptstyle #4}\ar[u]_{#6}}}

\newcommand{\syme}[8]
    {\xymatrix@C=20pt@R=15pt{
{\scriptstyle #1}\ar[d]_{#5}\ar[r]^{ #7}
&{\scriptstyle #2}
                      \\
{\scriptstyle #3}\ar[r]^{#8}&{\scriptstyle #4}\ar[u]_{#6}}}

\begin{document}
\title{Classification of linear mappings between indefinite inner product spaces}

\author[zai]{Juan Meleiro}
\ead{juan.meleiro@me.com}
\address[zai]{Instituto de Matem\'atica e Estat\'istica, Universidade de S\~ao Paulo, Brasil}

\author[ser]{Vladimir V.~Sergeichuk\corref{cor}}
\ead{sergeich@imath.kiev.ua}
\address[ser]{Institute of Mathematics,
Tereshchenkivska 3, Kiev, Ukraine}

\author[zai]{Thiago Solovera}
\ead{thiago.solovera.nery@usp.br}

\author[zai]{Andr\'{e} Zaidan}
\ead{andre.zaidan@gmail.com}

\cortext[cor]{Corresponding author.}


\begin{abstract}
Let $\mc A:U\to V$ be a linear mapping between vector spaces $U$ and $V$ over a field or skew field $\FF$ with symmetric, or skew-symmetric, or Hermitian forms
$
\mc B:U\times U\to\FF$ and $\mc C:V\times V\to\FF.$

We classify the triples $(\mc A,\mc B,\mc C)$ if $\FF$ is $\RR$, or $\CC$, or the skew field of quaternions $\HH$.
We also classify the triples $(\mc A,\mc B,\mc C)$ up to classification of symmetric forms and Hermitian forms if the characteristic of $\FF$ is not 2.
 \end{abstract}


\maketitle

\section{Introduction}
\label{s_intr}

We consider a triple
\begin{equation}\label{asa}
\mc A: U\to V,\qquad \mc B:U\times U\to\FF,\qquad\mc C:V\times V\to\FF
\end{equation}
consisting of
a linear mapping $\mc A$ and two forms
$\mc B$ and $\mc C$ on finite-dimensional vector spaces $U$ and $V$
over a field or skew field $\FF$ of characteristic not 2.
Each of the forms $\mc B$ and $\mc C$ is either symmetric or skew-symmetric if $\FF$ is a field, or both the forms are Hermitian with respect to a fixed nonidentity involution in $\FF$.

A canonical form of the triple of matrices of \eqref{asa} over a field $\FF$ of characteristic not 2 was obtained in the deposited manuscript \cite{ser_prep} up to classification of Hermitian forms over finite extensions of $\FF$. The aim of this paper is to give a detailed exposition of this result and extend it to triples \eqref{asa}  over a skew field of characteristic not 2.
We give canonical matrices of \eqref{asa} over $\RR$, $\CC$, and the skew field of quaternions $\HH$.

Other canonical matrices of \eqref{asa} with nonsingular forms $\mc B$ and $\mc C$ over the fields $\RR$ and $\CC$ were given by  Mehl, Mehrmann, and Xu \cite{meh,meh1,meh2}, and by Bolshakov and Reichstein \cite{bol}.

Following \cite{ser_prep}, we represent the triple \eqref{asa} by the graph
\begin{equation}\label{enr}
\begin{split}
\syk{ U}{V}{\mc A}{\mc B}{\mc C}{1}{2}{\varepsilon}{\delta }
\end{split}
\end{equation}
in which $\varepsilon =+$ if $\mc B$ is symmetric or Hermitian and  $\varepsilon =-$ if $\mc B$ is skew-symmetric; $\delta =+$ if $\mc C$ is symmetric or Hermitian and  $\delta =-$ if $\mc C$ is skew-symmetric.

Choosing bases in $U$ and $V$, we give \eqref{asa} by the triple $(A,B, C)$ of matrices
of $\mc A$, $\mc B$, and $\mc C$. Changing bases, we can reduce it by transformations
\begin{equation}\label{jyr}
(A,B,C)\mapsto (S^{-1}AR,R^{\as}BR,\,S^{\as}CS),\,
\end{equation}
in which $R$ and $S$ are nonsingular and
\begin{equation*}\label{wed}
M^{\as}=M^{\T}\quad\text{or}\quad M^{\as}=\widetilde{ M}^{\T}
\end{equation*}
with respect to a fixed
involution $a\mapsto \tilde a$ in $\FF$.
Thus, we consider the canonical form problem
for matrix triples under transformations \eqref{jyr}.
We represent the matrix triple $(A,B, C)$  by the graph
\begin{equation}\label{ehy}
\begin{split}
\syk{m}{n}{A}{B}{C}{1}{1}{\varepsilon }{\delta }\qquad
\xymatrix@R=13pt@C=45pt{
m:=\dim U,\\ n:=\dim V.}
\end{split}
\end{equation}

The \emph{direct sum} of matrices is $A\oplus B:=\matt{A&0\\0&B}$ and of matrix triples is
\begin{equation}\label{esw}
\begin{split}
\syk{m_1}{n_1}{A_1}{B_1}{C_1}{2}{4}{\varepsilon }{\delta }\
\syk{m_2}{n_2}{\hspace{-3mm}
{\textstyle \oplus}\ \
A_2}{B_2}{C_2}{2}{4}{\varepsilon }{\delta }\ \
\syk{m_1+m_2}{n_1+n_2}
{\hspace{-3mm}{\textstyle :=}\ \ A_1\oplus A_2}{B_1\oplus B_2}{C_1\oplus C_2}{5}{12}{\varepsilon }{\delta }
\end{split}
\end{equation}

The main result will be formulated in Section \ref{sab}. In the following theorem, we formulate it in the most important case: for linear mappings between indefinite inner product spaces (an \emph{indefinite inner product space} is a complex vector space with scalar product given by a nonsingular Hermitian form). Mehl, Mehrmann, and Xu \cite{meh} gave another classification of linear mappings between indefinite inner product spaces; their classification is presented in \cite[Section 6.5]{meh2}. This classification problem was also studied by Bolshakov and Reichstein \cite[Section 6]{bol}.
We refer the reader to Gohberg, Lancaster, and Rodman's book \cite{goh} for a recent account of the indefinite linear algebra.
\begin{theorem}\label{tt}
For each triple \eqref{enr} consisting of a
linear mapping $\mc A:U\to V$ and nonsingular
Hermitian forms $\mc B$ and $\mc C$ on complex
vector spaces $U$ and $V$,  there exist bases of
$U$ and $V$ in which the triple \eqref{ehy} of
matrices of $\mc A,\mc B,\mc C$ is a direct
sum, determined uniquely up to permutation of
summands,
of triples of the form
\begin{gather}\label{111}
\nonumber
\begin{split}
\syyk{r}{r}%
{I_r}{a  \zzz}%
{a  \matt{0&&&\lambda\\&&\lambda &1
\\[-1mm]&\iddots&\iddots\\\lambda
&1&&0}\ (0\ne\lambda \in\RR)}{1.5}{1.5}{+}{+}
          \qquad
\xymatrix@R=30pt@C=45pt{
{\scriptstyle 2r} \ar@{->}[d]_{\matt{I_r&0\\0&J_r(\mu)}}
^{\qquad\scriptstyle(\mu=\alpha +\beta i,\ \alpha ,\beta \in\RR,\ \beta >0)
}
 \save !<3 pt,0cm>
\ar@{-}@(ur,dr)^{\matt{0&I_r\\I_r&0}}_{+\,}
\restore
                      \\
{\scriptstyle 2r}
 \save !<3 pt,0cm>
\ar@{-}@(ur,dr)^{\matt{0&I_r\\I_r&0}}_{+\,}
\restore }
\end{split}
\\
\label{22a}
\begin{split}
\syyk{r}{r-1}
{\matt{1&0&&0\\[-1.5mm]&\ddots&\ddots
\\0&&1&0}\ \ }
{a  \zzz}
{a  \zzz}{1.5}{7}{+}{+}
      \qquad \qquad
\syyk{r-1}{r}
{\matt{\\1&&0
\\0&\smash{\ddots}
\\&\smash{\ddots}&1
\\
0&&0\\[1pt]}\ \ }
{a  \zzz}
{a  \zzz}{7}{1.5}{+}{+}
\end{split}
\end{gather}
in which
$r\in\{1,2,\dots\}$, $a  \in\{1,-1\}$, and $J_r(\mu)$ is an upper-triangular Jordan block.
\end{theorem}

Denote by $0_{pq}$ the $p\times q$ zero matrix with $p,q\in\{0,1,2,\dots\}$. Note that
\[
A_{mn}\oplus 0_{0q}=\mat{A_{mn}&0_{mq}\\0_{0n}&0_{0q}}
=\mat{A_{mn}&0_{mq}},\quad
A_{mn}\oplus 0_{p0}=\mat{A_{mn}&0_{m0}\\0_{pn}&0_{p0}}
=\mat{A_{mn}\\0_{pn}}
\]
for any $m\times n$ matrix $A_{mn}$.
The triples \eqref{22a} with $r=1$ have the form
\[
\syk{1}{0}{0_{01}}{[a  ]}{0_{00}}{1}{2}{+}{+}\qquad\qquad
\syk{0}{1}{0_{10}}{0_{00}}{[a  ]}{1}{2}{+}{+}
\]

We obtain our canonical form using the procedure  developed by Roiter and Sergeichuk in \cite{roi,ser_first,ser_prep,ser_izv} and presented in \cite{hor-ser_mixed,ser_brazil,ser_isom}.
If $\FF$ is a skew field (which can be commutative) of characteristic not 2 that is finite dimensional over its center, then this procedure reduces the problem of classifying any system of linear mappings and bilinear or sesquilinear forms over $\FF$ to the problems  of classifying (i) some system of linear mappings over $\FF$ and (ii) Hermitian forms over finite extensions of the center of $\FF$. The solution of problem (ii) is given by the law of inertia if $\FF$ is $\RR$, or $\CC$, or the skew field of quaternions $\HH$.

Over a field $\FF$ of characteristic not 2, Sergeichuk \cite{ser_prep}\footnote{This is a deposited manuscript. There were very few mathematical journals in the USSR and publications  abroad were forbidden. However,
a manuscript could be deposited in an  institute of information by recommendation of
a scientific council, which was considered as a publication. The abstracts of all deposited manuscripts on mathematics were published in the RZhMat, which is the Russian analogue of MathRev. Everyone can order a copy of each deposited manuscript.} obtained canonical forms,  up to classification of Hermitian forms over finite extensions of $\FF$, for matrices of
\begin{itemize}
  \item[(a)] bilinear and sesquilinear forms,

  \item[(b)] pairs of symmetric, or skew-symmetric, or Hermitian forms,

\item[(c)] self-adjoint or isometric operators in a space with scalar product given by a nonsingular form that is symmetric, or skew-symmetric, or Hermitian,

\item[(d)]
 systems represented by the graphs
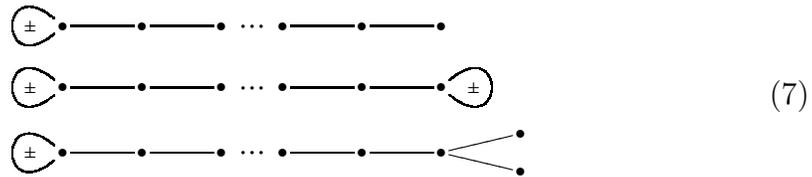
\begin{equation}\label{piug}
\begin{split}
&\xymatrix{*{\ci}\ar@{-}[r]
\ar@(ul,dl)@{-}^{\,
\pm}& *{\ci}
\ar@{-}[r]&*{{\ci}\
\cdots\ {\ci}}
\ar@{-}[r]&*{\ci}\ar@{-}[r]&
*{\ci} }
                                  \\[2mm]
&\xymatrix{*{\ci}\ar@{-}[r]
\ar@(ul,dl)@{-}^{\,
\pm}& *{\ci}
\ar@{-}[r]&*{{\ci}\
\cdots\ {\ci}}
\ar@{-}[r]&*{\ci}\ar@{-}[r]&
*{\ci}
\ar@(ur,dr)@{-}_{\pm\,}
}
                                     \\
&\xymatrix@R=0,5pt{
&&&&&*{\ci}\ar@{-}[dl]
\\*{\ci}\ar@{-}[r]
\ar@(ul,dl)@{-}^{\,
\pm}& *{\ci}
\ar@{-}[r]&*{{\ci}\
\cdots\ {\ci}}
\ar@{-}[r]&*{\ci}\ar@{-}[r]&
*{\ci}
\\
&&&&&*{\ci}\ar@{-}[ul]}\qquad\qquad
\end{split}
\end{equation}
in which each line between two vertices is an arrow $\longrightarrow$ or  $\longleftarrow$. The vertices are vector spaces over $\FF$, the arrows are linear mappings, and the loops are symmetric, or skew-symmetric, or Hermitian forms.
\end{itemize}

These canonical matrices were also published in \cite{ser_izv} for (a), (b), and (c), and in \cite{ser_rez} for (d). The procedure and the graphs \eqref{piug} are presented in the survey article \cite[Theorem 3.2]{ser_brazil};  see also \cite{hor-ser_transpose, hor-ser_can, hor-ser_bilin,ser_isom}.

Using the canonical matrices of the system
\raisebox{7pt}{$\xymatrix@R=0pt@C=8pt{
&*{\ci}\ar@{-}[dl]
\\*{\ci}
\ar@(ul,dl)@{-}^{\,
\pm}
\\
&*{\ci}\ar@{-}[ul]}
$}, Sergeichuk \cite{ser_sub} classified pairs of subspaces in a space with scalar product given by a symmetric, or skew-symmetric, or Hermitian form. The canonical matrices of the system
\ $\xymatrix@C=17pt{*{\ci}\ar@{->}[r]
\ar@(ul,dl)@{-}^{\,
\pm}& *{\ci}
\ar@(ur,dr)@{-}_{\pm\,}
}$ \ are given in the next section.

The procedure developed in \cite{roi,ser_prep} is based on Roiter's quivers with involution \cite{roi}, which were also studied by Derksen and Weyman \cite{der+wey} (they use the term ``symmetric quivers''), Bocklandt \cite{boc}, Shmelkin \cite{shm}, and Zubkov \cite{zub}.

\section{Main result}\label{sab}

Let $\FF$ be a skew field (which can be commutative) of characteristic not $2$ with a fixed \emph{involution} $a\mapsto \tilde a$; that is,
a bijection $\FF\to\FF$ satisfying
\[
\widetilde{a+b}=\tilde a+\tilde b,\quad
\widetilde{\,ab\,}\!=\tilde b\tilde a,\quad\Tilde{\Tilde a}=a
\]
for all $a,b\in\FF$. This involution can be the identity if $\FF$ is a field. Elementary linear algebra over a skew field can be found in
Bourbaki \cite[Chapter II]{bour}.

All vectors spaces that we consider are finite dimensional \emph{right} vector spaces over $\FF$. Each linear mapping $\mc A: U\to V$ satisfies $\mc A(ua+vb)=(\mc Au)a+(\mc Av)b$ and each sesquilinear form
 $\mc B:U\times U\to\FF$  satisfies
\[
\mc B(ua+vb,w)=\tilde a\mc B(u,w)+\tilde b\mc B(v,w),\quad
\mc B(u,va+wb)=\mc B(u,v)a+\mc B(u,w)b
\]
for all $a,b\in\FF$ and $u,v,w\in U$. A form $\mc B:U\times U\to\FF$ is \emph{Hermitian}  if $\mc B(u,v)=\widetilde {\mc B(v,u)}$ and \emph{skew-Hermitian} if $B(u,v)=-\widetilde {\mc B(v,u)}$ for all $u,v\in U$.

We consider a triple
\begin{equation}\label{as1a}
\begin{split}
\xymatrix@R=20pt@C=45pt{
{\scriptstyle U}\ar@{->}[d]_{\smash{\textstyle\mc T:}\qquad\mc A}^{\qquad\qquad\qquad
\textstyle\varepsilon,\delta \in\{+,-\}}
 \save !<1 mm,0cm>
\ar@{-}@(ur,dr)^{\mc B}_{\varepsilon\,}
\restore
                      \\
{\scriptstyle V}
 \save !<2 pt,0cm>
\ar@{-}@(ur,dr)^{\mc C}_{\delta\,}
\restore }
\end{split}
\end{equation}
that consists of a linear mapping $\mc A: U\to V$ between vector spaces over $\FF$ and two sesquilinear forms $\mc B:U\times U\to\FF$ and $\mc C:V\times V\to\FF$ satisfying
\begin{equation}\label{hsr}
\mc B(u,u')=\varepsilon \widetilde {\mc B(u',u)},\qquad
\mc C(v,v')=\delta \widetilde {\mc C(v',v)}
\end{equation}
for all $u,u',v,v'$.
For simplicity, \emph{we suppose that the forms $\mc B$ and $\mc C$ are Hermitian $($i.e., $\varepsilon=\delta=+)$  if the involution $a\mapsto \tilde a$ is not the identity.} This condition is not restrictive if $\FF$ is a field since then each skew-Hermitian form can be made Hermitian by multiplying by any $a-\tilde a\ne 0$; $a\in\FF$.

We say that two triples $\mc T$ and $\mc T'$ of the form \eqref{as1a} with the same $\varepsilon $ and $\delta $ are \emph{isomorphic} and write $\mc T\simeq\mc T'$ if there exist linear bijections $\varphi:U\to U'$ and $\psi:V\to V'$ that transform $\mc T$ to $\mc T'$:
\[
    \xymatrix@R=20pt@C=45pt{
{\scriptstyle U}
            \ar@/^1.2pc/@{-->}[r]^{\varphi}
\ar@{->}[d]_{\mc A}
 \save !<3 pt,0cm>
\ar@{-}@(ur,dr)^{\mc B}_{\varepsilon\,}
\restore
     & 
{\scriptstyle U'}\ar@{->}[d]^{\mc A'=\psi\mc A\varphi ^{-1}}
 \save !<4 pt,0cm>
\ar@{-}@(ur,dr)_{\varepsilon\, }^{\mc B'\qquad\textstyle \mc B(u_1,u_2)=\mc B'(\varphi u_1,\varphi u_2)\text{ for all }u_1,u_2\in U,}
\restore
                      \\
{\scriptstyle V}
             \ar@/_1.2pc/@{-->}[r]_{\psi}
 \save !<3 pt,0cm>
\ar@{-}@(ur,dr)^{\mc C}_{\delta\, }
\restore
          &  
{\scriptstyle V'}
 \save !<4 pt,0cm>
\ar@{-}@(ur,dr)_{\delta \,}^{\mc C'\qquad\textstyle \mc C(v_1,v_2)=\mc C'(\psi v_1,\psi v_2)\text{ for all }v_1,v_2\in V.}
\restore }
\]
The \emph{direct sum} of triples is defined as in \eqref{esw}:
\begin{equation*}\label{esw1}
\begin{split}
\syk{U_1}{V_1}{\mc A_1}{\mc B_1}{\mc C_1}{2}{4}{\varepsilon }{\delta }\
\syk{U_2}{V_2}{\hspace{-3mm}
{\textstyle \oplus}\ \
\mc A_2}{\mc B_2}{\mc C_2}{2}{4}{\varepsilon }{\delta }\
\syk{U_1\oplus U_2}{V_1\oplus V_2}{\hspace{-3mm}{\textstyle :=}\ \ \mc A_1\oplus\mc A_2}{\mc B_1\oplus \mc B_2}{\mc C_1\oplus\mc C_2}{5}{12}{\varepsilon }{\delta }
\end{split}
\end{equation*}

We say that a triple \eqref{as1a} is \emph{regular} if $\mc A$ is bijective. Each regular triple is isomorphic to a triple of the form
\begin{equation*}\label{enr1}
\begin{split}
\syk{ U}{U}{1}{\mc B}{\mc C}{1}{2}{\varepsilon}{\delta }
\end{split}
\end{equation*}
which we call \emph{strictly regular}.
We say that a triple is \emph{strictly singular} if it is not isomorphic to a direct sum with a regular direct summand.

For
$r\in\{0,1,2,\dots\}$, define the $r\times r$ matrices
\begin{equation}\label{asdq}
Z_r= Z_{r,+}:=\mat{0&&1\\&\iddots\\1&&0},
\quad\text{and}\
Z_{r,-}:=\mat{0&-Z_{r/2}\\Z_{r/2}&0}
\text{ if $r$ is even}.
\end{equation}
For
$r\in\{1,2,\dots\}$, define the $(r-1)\times r$ matrices
\begin{equation}\label{1.4a}
F_r:=\mat{1&0&&0\\&\ddots&\ddots&\\0&&1&0},\qquad
G_r:=\mat{0&1&&0\\&\ddots&\ddots&\\0&&0&1}.
\end{equation}
In particular, $F_1=G_1=0_{01}$.

Let
\begin{equation}\label{fdf}
\xymatrix{\ar@{-}@(ur,dr)
^{A^{\oslash}}_{\varepsilon \,}}\quad\text{denote}\quad
\xymatrix{\ar@{-}@(ur,dr)
^{\matt{0&\varepsilon A^{\as}\\A&0}}_{\varepsilon\,}}
\end{equation}
in which $\varepsilon \in\{+,-\}$ and $A^{\as}=\widetilde A^{\T}$ is the \emph{adjoint} matrix.

We can now formulate our main result, which was given over a field in \cite[p.\,44]{ser_prep}.

\begin{theorem}\label{tttt}
Let $\FF$ be a skew field $($which can be commutative$)$ of characteristic not $2$ with a fixed involution $a\mapsto \tilde a$. Let us fix $\varepsilon,\delta \in\{+,-\}$ if the involution is the identity, and put $\varepsilon=\delta= +$ if the involution is not the identity.

{\rm(a)} Let
\eqref{as1a} be a triple,
consisting of a linear mapping $\mc A: U\to V$ between right vector spaces over $\FF$ and two sesquilinear forms $\mc B:U\times U\to\FF$ and $\mc C:V\times V\to\FF$ satisfying \eqref{hsr}. Then the triple \eqref{as1a}
is isomorphic to a direct sum of a strictly regular triple and a strictly singular triple; these summands are uniquely determined, up to isomorphism.

{\rm(b)} Two strictly regular triples are isomorphic if and only if their forms are simultaneously equivalent:
\begin{equation*}\label{enr3}
\begin{split}
\syk{ U}{U}{1}{\mc B}{\mc C}{1}{2}{\varepsilon}{\delta }\
               {\xymatrix@=6pt{{}\\ \simeq}} \
\syk{ U'}{U'}{1}{\mc B'}{\mc C'}{1.5}{2.5}{\varepsilon}{\delta }
\quad         {\xymatrix@=6pt{{}\\ \Longleftrightarrow}}\ \quad
\xymatrix@R=5pt{\\
{\scriptstyle U}
 \save !<1.5pt,0cm>
\ar@{-}@(ur,dr)_{\delta\,}^{\mc C}
\restore
 \save !<-2.5pt,0cm>
\ar@{-}@(ul,dl)^{\;\varepsilon}_{\mc B}
\restore
 }
\  {\xymatrix@=6pt{{}\\ \simeq}} \
\xymatrix@R=6pt{\\
{\scriptstyle U'}
 \save !<2.5pt,0cm>
\ar@{-}@(ur,dr)_{\delta\,}^{\mc C'}
\restore
 \save !<-3.5pt,0cm>
\ar@{-}@(ul,dl)^{\;\varepsilon}_{\mc B'}
\restore
 }
\end{split}
\end{equation*}

{\rm(c)}
Each  strictly singular triple \eqref{as1a} possesses bases of\/
$U$ and $V$, in which the triple \eqref{ehy} of its
matrices is a direct
sum
of triples of the types
\begin{equation}\label{b1a}
\begin{split}
\syrk{r}{r-1}
{\smash{{\displaystyle F_r^1(a):}\ }F_{r}}
{aZ_{r,\varepsilon}}
{aZ_{r-1,\delta}}
{0}{4}
{\qquad \begin{smallmatrix}
(\text{$r$ is even if $\varepsilon =-$,\ }\\
\text{$r$ is odd if $\delta =-$})
        \end{smallmatrix}
}
                               \qquad
\syrk{2r}{2r-2}
{F_{r}\oplus G_r}
{I_{r}\od}
{I_{r-1}\od}{1}{5}
{\qquad \begin{smallmatrix}
(\text{$r$ is odd and $\varepsilon =-$,\quad\ \ }\\
\text{\, or $r$ is even and $\delta  =-$})
        \end{smallmatrix}}
\end{split}
\end{equation}
\begin{equation}\label{b2}
\begin{split}
\syrk{r-1}{r}
{\smash{{\displaystyle F_r^{\T}(a):}\  }F_{r}^{\T}}{aZ_{r-1,\varepsilon }}
{aZ_{r,\delta }}{2}{0}
{\qquad \begin{smallmatrix}(\text{$r$ is odd if $\varepsilon =-$,\quad}\\ \
\text{$r$ is even if $\delta =-$})\
        \end{smallmatrix}
}
                             \qquad
\syrk{2r-2}{2r}
{F_{r}^{\T}\oplus G_r^{\T}}
{I_{r-1}\od
}
{I_{r}\od
}{2}{2}
{\qquad \begin{smallmatrix}
(\text{$r$ is even and $\varepsilon =-$,\quad}\\
\text{ or $r$ is odd and $\delta  =-$})
        \end{smallmatrix}
}
\end{split}
\end{equation}
\begin{equation}\label{b5}
\begin{split}
\syk{2r}{2r-1}{I_r\oplus F_r}{I_r\od}{G_r\od}{1}{8}{\varepsilon}{\delta}
        \hspace{85pt}
\syk{2r-1}{2r}{I_r\oplus F_r^{\T}}{G_r\od}{I_r\od}{3}{2.5}
{\varepsilon}{\delta}
\end{split}
\end{equation}
\begin{equation}\label{b6}
\begin{split}
\syk{2r-1}{2r-2}{I_{r-1}\oplus F_r}{(G_r^{\T})\od}{I_{r-1}\od}{3}{8}
{\varepsilon}{\delta}
         \hspace{25pt}
\syk{2r-2}{2r-1}{I_{r-1}\oplus F_r^{\T}}{I_{r-1}\od}{(G_r^{\T})\od}{3}{8}
{\varepsilon}{\delta}
         \hspace{25pt}
\syk{2r}{2r}{I_{r}\oplus J_r(0)}{I_{r}\od}{I_r\od}{1}{2.5}
{\varepsilon}{\delta}
\end{split}
\end{equation}
in which
$r\in\{1,2,\dots\}$  and $0\ne a=\tilde a\in\FF$.

The summands of types $F_r^{\sigma}(a)$ with $\sigma \in\{1,\T\}$  are
determined up to replacement of
the whole group of summands
\begin{equation*}\label{hge}
F_r^{\sigma}(a_1)\oplus\dots\oplus F_r^{\sigma}(a_k)
\end{equation*}
with the same $r$ and $\sigma $
by any direct sum
\[
F_r^{\sigma}(b_1)\oplus\dots\oplus F_r^{\sigma}(b_k),\qquad b_1=\tilde b_1,\dots,b_k=\tilde b_k\in\FF\smallsetminus\{0\}
\]
with the same $r$ and $\sigma $
such that the
Hermitian forms
\begin{equation*}\label{777}
\tilde x_1a_1x_1+\dots+
\tilde x_ka_kx_k,
\qquad
\tilde x_1b_1x_1+\dots+
\tilde x_kb_kx_k
\end{equation*}
are equivalent over
$\FF$. The other summands are uniquely determined up to permutation.

{\rm (d)} Let $\FF$ be $\RR$, or $\CC$, or the skew field of quaternions $\HH$. Then each  strictly singular triple \eqref{as1a} possesses bases of\/ $U$ and $V$, in which the triple \eqref{ehy} of its
matrices is a direct
sum, uniquely determined up to permutation of summands, of triples of types
\eqref{b1a}--\eqref{b6}, in which
\begin{itemize}
  \item $a=1$ if $\FF$ is $\CC$ with the identity involution, or $\HH$ with involution that differs from the quaternion conjugation
\begin{equation}\label{tgh}
\alpha +\beta i+\gamma j+\delta k\ \mapsto\ \alpha -\beta i-\gamma j-\delta k\qquad(\alpha,\beta,\gamma,\delta\in\RR);
\end{equation}

  \item $a\in\{-1,1\}$ if $\FF$ is $\RR$, or $\CC$ with complex conjugation, or $\HH$ with the quaternion conjugation.
\end{itemize}
\end{theorem}

The triples \eqref{b1a}--\eqref{b6} are the triples $(A_2''1)$, $(A_2''2)$, $(A_2''2)_a$, and $(A_2''3)$ from \cite[p.\,44]{ser_prep} and the  triples that are ``dual'' to them (they are obtained by replacing the vector spaces by the dual vector spaces  \eqref{aan} and the linear mappings by the adjoint mappings \eqref{qmi}; the triples with forms on the dual spaces are also classified in \cite{ser_prep}).
Each involution in $\HH$ is either
\eqref{tgh}, or $\alpha +\beta i+\gamma j+\delta k\; \mapsto\; \alpha +\beta i+\gamma j-\delta k$ in a suitable set $\{i,j,k\}$ of orthogonal
imaginary units; see \cite[Theorem 2.4.4(c)]{rod_quat_book}.

Thus, Theorem \ref{tttt}
reduces the problem of classifying triples \eqref{as1a}
over $\FF$ up to isomorphism
\begin{itemize}
  \item[(i)] to the problem of classifying pairs of forms $\xymatrix@R=6pt{*{\ci}
\ar@{-}@(ur,dr)_{\pm\,}
\ar@{-}@(ul,dl)^{\;\pm}
}$ over $\FF$,
 and
  \item[(ii)] (if $\FF$ is not $\RR$, $\CC$, and $\HH$) to the problem of classifying Hermitian forms over $\FF$.
\end{itemize}
We do not consider the problem (i); its solution is given  in \cite[Theorem 4]{ser_izv} over a field $\FF$ of characteristic not $2$ up to classification of Hermitian forms over finite extensions of $\FF$, which gives its full solution if $\FF$ is $\RR$ or $\CC$ due to the law of inertia. The pairs $\xymatrix@R=6pt{*{\ci}
\ar@{-}@(ur,dr)_{\pm\,}
\ar@{-}@(ul,dl)^{\;\pm}
}$ over $\RR$ and $\CC$ are also classified in \cite{hor-ser_can,lan1,lan2,thom} and other papers. The pairs $\xymatrix@R=6pt{*{\ci}
\ar@{-}@(ur,dr)_{\pm\,}
\ar@{-}@(ul,dl)^{\;\pm}
}$ over $\HH$ are classified in \cite{kar,rodm1,rodm2,rod_quat_book}.

\section{The reducing procedure}\label{red}

The procedure that reduces the problem of classifying systems of linear mappings and forms to the problem of classifying systems of linear mappings is described in \cite[Section 1]{ser_izv} and \cite[Section 3.5]{ser_isom}.
In this section, we present it for the problem of classifying triples \eqref{as1a}.

For each right vector space $V$ over $\FF$,
  \begin{equation}\label{aan}
\begin{matrix}
\text{$V^{\as}$ is the right \emph{dual space} of \emph{semilinear functionals}},\\
\text{that is, mappings $\varphi :V\to \FF$ satisfying}\\
\text{$\varphi (ua+vb)=\tilde a\varphi (u)+\tilde b\varphi (v)$ for all $a,b\in\FF$ and $u,v\in V$.}
\end{matrix}
\end{equation}
For each linear mapping $\mc A:U\to V$,
\begin{equation}\label{qmi}
\begin{matrix}
\text{$\mc A^{\as}:V^{\as}\to U^{\as}$ is the \emph{adjoint mapping} defined by} \\
(\mc A^{\as}\varphi )(u):=\varphi (\mc A u)\text{ for all $u\in U$ and $\varphi \in V^{\as}$.}
\end{matrix}
\end{equation}

Each sesquilinear form $\mc B:U\times V\to \FF$ defines
\begin{itemize}
  \item
the linear mapping (we denote it by the same letter)
\begin{equation}\label{nnx}
\mc B: V\to U^{\as},\qquad v\mapsto \mc B(?,v),
\end{equation}
  \item
the adjoint linear mapping
\[\mc B^{\as}:U\to V^{\as},\qquad u\mapsto \widetilde{\mc B(u,?)}.\]
  \end{itemize}
If the form $\mc B$ is Hermitian, then the mapping \eqref{nnx} is self-adjoint (i.e., $\mc B=\mc B^{\as}$).

Thus, the triple $\mc T$ in \eqref{as1a} defines in a one-to-one manner the \emph{quadruple of linear mappings}
\begin{equation}\label{vri}
\begin{split}
\symaa{U}{U^{\as}}
{V}{V^{\as}}
{\smash{\textstyle\underline{\mc T}:}\qquad\mc A}{\mc A^{\as}}{\mc B=\varepsilon \mc B^{\as}}
{\mc C=\delta \mc C^{\as}}
\end{split}
\end{equation}
(in terms of \cite{roi,ser_izv, hor-ser_mixed}, the quadruple \eqref{vri} is a self-dual representation of the quiver
\begin{equation}\label{kjr}
\begin{split}
\symaa{u}{u^*}
{v}{v^*}
{\smash{\textstyle \underline G:}\qquad \alpha }{\alpha^*}{\beta=\varepsilon \beta^*}
{\gamma=\delta \gamma^*}
\end{split}
\end{equation}
with involutions in the set of vertices and in the set of arrows).

Thus, we can consider quadruples of the form \eqref{vri} instead of the triples \eqref{as1a}.
We will classify the quadruples \eqref{vri} using the classification of \emph{arbitrary} quadruples of linear mappings
\begin{equation}\label{vku}
\begin{split}
\symaa{U_1}{U_2}
{V_1}{V_2}
{\smash{\textstyle \mc P:}\qquad \mc A_1}{\mc A_2}{\mc B}
{\mc C}
\end{split}
\end{equation}
(that is, representations of the quiver \eqref{kjr}, which we consider as a quiver without involutions).

The vector
\begin{equation}\label{tno}
\dim P:=(\dim U_1,\dim U_2, \dim V_2,\dim V_1)
\end{equation}
is called the \emph{dimension} of \eqref{vku}.

A \emph{homomorphism}
\begin{equation}\label{gfh}
\phi=
       \left(\begin{smallmatrix}
         \varphi_1&\varphi_2\\ \psi_1&\psi_2
       \end{smallmatrix}\right)
     :\ \mc P\to \mc P'
\end{equation}
of quadruples $\mc P$ and $\mc P'$
is a sequence $\varphi_1,\varphi_2, \psi_1,\psi_2$ of linear mappings
\begin{equation*}\label{vjo3}
\begin{split}
    \xymatrix@R=15pt@C=25pt{
{\scriptstyle U_1}
       \ar@/^1.5pc/@{-->}[rrr]^{\varphi _1}
\ar[d]_{\smash{\textstyle\phi:}\quad \mc A_1}\ar[r]^{\mc B}
&{\scriptstyle U_2}
       \ar@/^1.5pc/@{-->}[rrr]^{\varphi _2}
&&
{\scriptstyle U_1'}
\ar[d]_{\mc A'_1}\ar[r]^{\mc B'\ }
&{\scriptstyle U'_2}
                                      \\
{\scriptstyle V_1}
             \ar@/_1.5pc/@{-->}[rrr]_{\psi_1}
\ar[r]^{\mc C}&{\scriptstyle V_2}
             \ar@/_1.5pc/@{-->}[rrr]_{\psi_2}
\ar[u]_{\mc A_2}&&
{\scriptstyle V'_1}\ar[r]^{\mc C'}&{\scriptstyle V'_2}\ar[u]_{\mc A'_2}
}
\end{split}
\end{equation*}
such that
\[
\psi _1\mc A_1=\mc A_1'\varphi _1,\quad
\varphi_2\mc A_2=\mc A_2'\psi _2,\quad
\varphi_2\mc B=\mc B'\varphi _1,\quad
\psi_2\mc C=\mc C'\psi_1.
\]
A homomorphism \eqref{gfh} is called an \emph{isomorphism} (we write $\phi: \mc P\:\mathaccent\sim\to\: \mc P'$ or $\mc P\simeq \mc P'$) if $\varphi _1,\varphi _2,\psi _1,\psi _2$ are bijections.

Let $\varepsilon $ and $\delta $ be as in \eqref{as1a}.
For each quadruple $\mc P$ of the form \eqref{vku}, we define the \emph{dual quadruple}
\begin{equation}\label{vsw}
\begin{split}
\symaa{U_2^{\as}}{U_1^{\as}}
{V_2^{\as}}{V_1^{\as}}
{\smash{\textstyle \mc P^{\circ}:}\qquad \mc A_2^{\as}}{\mc A_1^{\as}}{\varepsilon \mc B^{\as}}{\delta \mc C^{\as}}
\end{split}
\end{equation}
The quadruple \eqref{vri} is self-dual.

For each homomorphism \eqref{gfh}, we define
the \emph{dual homomorphism} $\phi^{\circ}:=\left(
\begin{smallmatrix}\varphi_2^{\as} &\varphi_1^{\as}\\ \psi_2^{\as}&\psi_1^{\as}
\end{smallmatrix}\right): \mc P^{\prime\circ}\to \mc P^{\circ}$ between the dual quadruples:
\begin{equation*}\label{vjo4}
\begin{split}
    \xymatrix@R=15pt@C=25pt{
{\scriptstyle U_2^{\prime\as}}
       \ar@/^2pc/@{-->}[rrr]^{\varphi _2^{\as}}
\ar[d]_{\smash{\textstyle\phi^{\circ}:}\quad \mc A_2^{\prime\as}}\ar[r]^{\varepsilon \mc B^{\prime\as}}
&{\scriptstyle U_1^{\prime\as}}
       \ar@/^2pc/@{-->}[rrr]^{\varphi _1^{\as}}
&&
{\scriptstyle U_2^{\as}}
\ar[d]_{\mc A_2^{\as}}\ar[r]^{\varepsilon \mc B^{\as}}
&{\scriptstyle U_1^{\as}}
                     \\
{\scriptstyle V_2^{\prime\as}}
\ar@/_1.5pc/@{-->}[rrr]_{\psi_2^{\as}}
\ar[r]^{\delta \mc C^{\prime\as}}&{\scriptstyle V_1^{\prime\as}}
\ar@/_1.5pc/@{-->}[rrr]_{\psi_1^{\as}}
\ar[u]_{\mc A_1^{\prime\as}}&&
{\scriptstyle V_2^{\as}}\ar[r]^{\delta \mc C^{\as}}&{\scriptstyle V_1^{\as}}\ar[u]_{\mc A_1^{\as}}
}
\end{split}
\end{equation*}

Define the \emph{direct sum} of quadruples:
\[
\symaa{U_1\oplus U_1'}{U_2\oplus U_2'}
{V_1\oplus V_1'}{V_2\oplus V_2'}
{\smash{\textstyle \mc P\oplus \mc P':}\qquad \mc A_1\oplus A_1'}{\mc A_2\oplus A_2'}{\mc B\oplus B'}
{\mc C\oplus C'}
\]

A quadruple is \emph{indecomposable} if it is not isomorphic to a direct sum of quadruples of smaller dimensions.
Let $\ind(\underline{G})$ be any set of indecomposable quadruples such that each indecomposable quadruple is isomorphic to exactly one quadruple from $\ind(\underline{G})$ (we give $\ind(\underline{G})$ in Lemma \ref{vbe2}). The procedure of constructing indecomposable canonical triples \eqref{as1a} consists of three steps; see details in \cite[\S 1]{ser_izv} and \cite[Section 3]{ser_isom}.

\begin{description}
  \item[Step 1.]
We replace each quadruple in $\ind (\underline{G})$ that is isomorphic to a self-dual quadruple by a self-dual quadruple. Let $\ind_0(\underline{G})\subset\ind(\underline{G})$ be the set of obtained self-dual quadruples $\mc P= \mc P^{\circ}$.
Denote by $\ind_1(\underline{G})$ a set consisting of
\begin{itemize}
                 \item
all $\mc Q\in\ind(\underline{G})\smallsetminus \ind_0(\underline{G})$ such that $\mc Q\simeq \mc Q^{\circ}$, and
                 \item
one quadruple from each pair $\{\mc Q,\mc R\}\subset\ind(\underline{G})$ such that $\mc Q\ne \mc R$ and $\mc Q^{\circ}\simeq \mc R$.           \end{itemize}

We have obtained a new set
$\ind (\underline{G})$
partitioned into 3 subsets:
\begin{equation}\label{4.8d}
{\ind (\underline{G})}
=
\begin{tabular}{|c|c|}
  \hline
\multicolumn{2}{|c|}
{${\cal P}={\cal P}^{\circ}\vphantom{{\hat{N}}}$}\\
 \hline &\\[-12pt]
  $\; {\cal Q}\;  $&
  ${\cal Q}^{\circ}\text{ if
  ${\cal Q}^{\circ}\not
  \simeq{\cal Q}$}$
  \\  \hline
\end{tabular}\,,\ \
\begin{matrix}
 {\cal P}\in
\ind_0(\underline{G}),\\[1pt]
{\cal
Q}\in\ind_1(\underline{G}).
\end{matrix}
\end{equation}

  \item[Step 2.] Let ${\cal P}\in
\ind_0(\underline{G})$.
Since $\mc P$ is an indecomposable quadruple, the algebra $\End(\mc P)$ of its endomorphisms is local, its radical $R$ consists of all non-invertible endomorphisms (see \cite[Lemma 1]{ser_izv}), and so $\TT(\mc P):=\End(\mc P)/R$ is a skew field. The mapping
\begin{equation*}\label{gmw}
\phi+R\mapsto(\phi+R)^{\circ}:=\phi^{\circ}+R,\qquad \phi\in\End(\mc P),
\end{equation*}
is an involution in $\TT(\mc P)$.
For each self-dual automorphism $\phi=\phi^{\circ}:=
\mattc{\varphi&\varphi^{\as}\\ \psi&\psi^{\as}}: \mc P\:\mathaccent\sim\to\: \mc P$, we denote by $\mc P^{\phi}$ the self-dual quadruple such that
\begin{equation}\label{ebh}
\mc P\
   \xrightarrow{\mattc{\varphi&1\\ \psi&1}}\
\mc P^{\phi}\
  \xrightarrow{\mattc{1&\varphi^{\as}\\ 1&\psi^{\as}}}\
\mc P.
\end{equation}
For each $0\ne a=a^{\circ}\in \TT(\mc P)$, we fix a  self-dual  automorphism $\phi_a=\phi_a^{\circ}\in a$ (we can take $\phi_a:=(\phi+\phi^{\circ})/2$ for any $\phi\in a$) and denote by $\mc P^a$ the triple \eqref{as1a} that corresponds to the self-dual quadruple $\mc P^{\phi_a}$.
For each Hermitian form
\[
f(x)=x^{\circ}_1a_1x_1+\dots+
x^{\circ}_ka_rx_k,\qquad
a_1=a_1^{\circ},\dots,a_k=a_k^{\circ}\in \TT(\mc P)\smallsetminus\{0\},
\]
over the skew field $\TT(\mc P)$, we put
\[
\mc P^{f(x)}:=
\mc P^{a_1}\oplus\cdots\oplus
\mc P^{a_k}.
\]

  \item[Step 3.]
For each quadruple $\mc Q\in\ind_1(\underline{G})$, we take the direct sum
\[
\symaaa{U_1\oplus  U_2^{\as}}{U_2\oplus U_1^{\as}}
{V_1\oplus  V_2^{\as}}{V_2\oplus  V_1^{\as}}
{{\textstyle \mc Q\oplus \mc Q^{\circ}:}\qquad\matt{\mc A_1&0\\0&\mc A_2^{\as}}}{\matt{\mc A_2&0\\0&\mc A_1^{\as}}}{\matt{\mc B&0\\0&\varepsilon \mc B^{\as}}}
{\matt{\mc C&0\\0&\mc \delta C^{\as}}}
\]
and make it self-dual by interchanging the summands in $U_2\oplus U_1^{\as}$ and in $V_2\oplus V_1^{\as}$:
\begin{equation*}\label{nmx}
\symaaa{U_1\oplus U_2^{\as}}{ U_1^{\as}\oplus U_2}
{V_1\oplus V_2^{\as}}{ V_1^{\as}\oplus V_2}
{\matt{\mc A_1 &0\\0&\mc A_2^{\as}}}{\matt{\mc A_1^{\as}&0\\0&\mc A_2}}{\matt{0&\varepsilon \mc B^{\as}\\\mc B&0}}
{\matt{0&\delta \mc C^{\as}\\\mc C&0}}
\end{equation*}
The corresponding triple is
\begin{equation}\label{ekx}
\begin{split}
\syk{U_1\oplus U_2^{\as}}{V_1\oplus V_2^{\as}}{\smash{\displaystyle\mc Q^+:}\qquad \matt{\mc A_1 &0\\0&\mc A_2^{\as}}}{\matt{0&\varepsilon \mc B^{\as}\\\mc B&0}}{\matt{0&\delta \mc C^{\as}\\\mc C&0}}{5}{14}{\varepsilon}{\delta }
                 \end{split}
\end{equation}
\end{description}

The following lemma is a special case of \cite[Theorem 1]{ser_izv} (see also \cite[Theoren 3.1]{ser_isom}) about \emph{arbitrary} systems of linear mappings and forms.

\begin{lemma}\label{jux}
Over a skew field
\/$\FF$ of characteristic not
$2$, each triple \eqref{as1a}
is isomorphic to a direct sum
\begin{equation}\label{ser03}
\mc P_1^{f_1(x)}\oplus\dots\oplus
\mc P_m^{f_m(x)}\oplus
\mc Q_1^+\oplus\dots\oplus\mc Q_n^+,
\end{equation}
in which $\mc P_1,\dots, \mc P_m\in\ind_0({\underline
G})$, $\mc P_i\ne\mc P_{i'}$ if
$i\ne i'$, each $f_i(x)$ is a Hermitian form over the skew field $\TT(\mc P_i)$, and $\mc Q_1,\dots,\mc Q_n\in\ind_1({\underline
G})$. This sum is
uniquely determined, up to
permutation of summands and
replacement of each
$\mc P_i^{f_i(x)}$ by
$\mc P_i^{g_i(x)}$, where
$f_i(x)$ and
$g_i(x)$ are equivalent
Hermitian forms over $\TT(\mc P_i)$.
\end{lemma}

If a skew field $\FF$ is finite-dimensional over its center $Z$ and $\mc P$ is a self-dual indecomposable quadruple, then $\TT(\mc P)$ is finite-dimensional over $Z$ under the natural imbedding of $Z$ into the center of $\TT(\mc P)$ and
the involution in $\TT(\mc P)$ extends the involution $a\mapsto \tilde a$ in $Z$.

If we have chosen a basis $e_1,\dots,e_n$ in a vector space $V$, then we always choose in the  dual space $V^{\as}$ the \emph{dual basis} consisting of the semilinear functionals $e_1^*,\dots,e_n^*:V\to\FF$ such that $e_i^*(e_j)=0$ if $i\ne j$ and $e_i^*(e_i)=1$ for all $i$.
\begin{equation}\label{njg}
\parbox[c]{0.8\textwidth}{If $A$ is the matrix of a linear mapping $\mc A:U\to V$ in some bases of $U$ and $V$, then $A^{\as}:=\widetilde A^{\T}$ is the matrix of the adjoint mapping $\mc A^{\as}:V^{\as}\to U^{\as}$  (see \eqref{qmi}) in the dual bases.
If $\mc B:U\times V\to \FF$ is a sesquilinear form, then its matrix in any bases of $U$ and $V$ coincides with the matrix of $\mc B:V\to U^{\as}$, $v\mapsto \mc B(?,v)$, in the same basis of $V$ and the dual basis of $U^{\as}$.
}
\end{equation}
Thus,
\begin{equation}\label{ekx1}
\begin{split}
\sykl{U_1\oplus U_2^{\as}}{V_1\oplus V_2^{\as}}{\smash{\displaystyle\mc Q^+:}\ \ \matt{\mc A_1 &0\\0&\mc A_2^{\as}}}{\matt{0&\varepsilon \mc B^{\as}\\\mc B&0}}{\matt{0&\delta \mc C^{\as}\\\mc C&0}}{5}{14}{\varepsilon}{\delta }
\sykl{m_1+m_2}{n_1+n_2}
{\smash{\displaystyle
\text{corresponds to } Q^+:}\ \
\matt{A_1 &0\\0& A_2^{\as}}}
{B\od}{C\od}
{5}{12}{\varepsilon}{\delta}
 \end{split}
\end{equation}
(see \eqref{fdf}), in which $m_1,m_2,n_1,n_2$ are the dimensions of $U_1,U_2,V_1,V_2$.

\section{Classification of quadruples of linear mappings}

Choosing bases in the spaces of a quadruple \eqref{vku}, we can give it by the quadruple of its matrices
\begin{equation}\label{vrj}
\begin{split}
\symaa{\FF^{m_1}}{\FF^{m_2}}
{\FF^{n_1}}{\FF^{n_2}}
{\smash{\textstyle P:}\qquad A_1}{A_2}{B}{C}
\end{split}
\end{equation}
This quadruple is isomorphic to the quadruple \eqref{vku}.
For abbreviation, we usually omit ``$\FF$'' in \eqref{vrj} (as in \eqref{ehy}). Two matrix quadruples are isomorphic (which means that
 \begin{equation*}\label{vre}
\begin{split}
\symaa{m_1}{m_2}
{n_1}{n_2}
{\smash{\textstyle\left(
\begin{smallmatrix}
\Phi_1 & \Phi_2 \\
\Psi_1& \Psi_2
\end{smallmatrix}\right)
:}\quad
A_1}{A_2}{B}{C}
\symaa{m_1}{m_2}
{n_1}{n_2}
{\ {\textstyle \mathaccent\sim\to}\ \
\Psi_1A_1\Phi_1^{-1}}{\Phi_2A_2\Psi_2^{-1}}
{\Phi_2B\Phi_1^{-1}}{\Psi_2C\Psi_1^{-1}}
\end{split}
\end{equation*}
for some nonsingular matrices $\Phi_1,\Phi_2,\Psi_1,\Psi_2$)
if and only if they give the same quadruple \eqref{vku} in different bases.

By \eqref{vsw} and \eqref{njg}, the \emph{dual quadruple} to \eqref{vrj} is the quadruple
\begin{equation*}\label{vswb}
\begin{split}
\symaa{m_2}{m_1}
{n_2}{n_1}
{\smash{\textstyle P^{\circ}:}\qquad A_2^{\as}}{A_1^{\as}}{B^{\as}}{ C^{\as}}
\end{split}
\end{equation*}
in which $A^{\as}:=\widetilde A^{\T}$ is the adjoint matrix.

A quadruple \eqref{vku} is a special case of a cycle of linear mappings
\begin{equation}\label{jsttw}
\begin{split}
\xymatrix{
{V_1}
\ar@{-}@/_1.5pc/[rrrr]^{\mathcal A_t}
\ar@{-}[r]^{\mathcal A_1}&
V_2\ar@{-}[r]^{\mathcal A_2\ } &
{\ \cdots\ }&{V_{t-1}}
\ar@{-}[l]_{\mathcal A_{t-2}}
\ar@{-}[r]^{\ \mathcal A_{t-1}}&{V_t}}
\end{split}
\end{equation}
in which each line is
$\longrightarrow$ or $\longleftarrow$, $V_1,\dots,V_t$ are vector spaces, and $\mc A_1,\dots,\mc A_t$ are linear mappings.
A canonical form of matrices of a cycle over a field is well known;
see, for example, \cite[Theorem 3.2]{ser-cycl}.
The regularizing algorithm from \cite{ser-cycl}
constructs for each cycle \eqref{jsttw}  its decomposition into a direct sum, in which the first summand is of the form
\begin{equation}\label{ttw}
\begin{split}
\xymatrix{
{V}
\ar@{-}@/_1.5pc/[rrrr]^{\mathcal A}
\ar@{-}[r]^{1}&
V\ar@{-}[r]^{1\ } &
{\ \cdots\ }&{V}
\ar@{-}[l]_{\ 1}
\ar@{-}[r]^{1}&{V}}
\end{split}
\end{equation}
with a bijective $\mc A$
and each other summand is an indecomposable canonical cycle that contains a nonbijective linear mapping. The proof of this algorithm is given over a field but it holds over the skew field $\FF$ too. Thus,
a canonical form of the cycle \eqref{ttw} over $\FF$ is obtained from the canonical form of a linear operator over a skew field, which is given in \cite[Chapter 3, Section 12]{jac}.
This ensures the following lemma.

\begin{lemma}\label{vbe2}
Let $\FF$ be a skew field, which can be commutative.
For each quadruple
\eqref{vku} over\/ $\FF$, there exist bases of $U_1,U_2,V_1,V_2$, in which the quadruple \eqref{vrj} of its matrices is a direct sum, uniquely determined up to permutations of summands, of indecomposable quadruples of the following types:
\begin{equation}\label{las1}
\begin{split}
\syme{r}{r}{r-1}{r-1}{F_r}{G_r^{\T}}
{I_r}{I_{r-1}}\qquad\qquad
\syme{r-1}{r-1}{r}{r}{F^{\T}_r}{G_r}{I_{r-1}}{I_r}
\end{split}
\end{equation}
\begin{equation}\label{las2}
\begin{split}
\syme{r}{r}{r}{r-1}{I_r}
{F^{\T}_r}{I_r}{G_r}\quad
\syme{r}{r-1}{r}{r}{I_r}{F_r}{G_r}{I_r}\quad
\syma{r-1}{r}{r-1}{r-1}{I_{r-1}}
{F^{\T}_r}{G^{\T}_r}{I_{r-1}}\quad
\syma{r-1}{r-1}{r-1}{r}{I_{r-1}}
{F_r}{I_{r-1}}{G^{\T}_r}
\end{split}
\end{equation}
\begin{equation}\label{las}
\begin{split}
\syme{r}{r}{r}{r}{I_r}{\Phi_r}{I_r}{I_r}\quad
\syme{r}{r}{r}{r}{J_r(0)}{I_r}{I_r}{I_r}\quad
\syme{r}{r}{r}{r}{I_r}{I_r}{J_r(0)}{I_r}\quad
\syme{r}{r}{r}{r}{I_r}{I_r}{I_r}{J_r(0)}
\end{split}
\end{equation}
and the quadruples dual to \eqref{las2}, in which $r\in\{1,2,\dots\}$, $F_r$ and $
G_r$ are defined in \eqref{1.4a}, and $\Phi_r$ is an $r\times r$ indecomposable canonical matrix under similarity over $\FF$.
\end{lemma}

Note that the dimensions \eqref{tno} of indecomposable quadruples from Lemma \ref{vbe2} are
\begin{equation*}\label{hgf}
(r,r,r,r),\ (r,r,r,r-1),\ (r,r,r-1,r-1),\ (r,r-1,r-1,r-1)
\end{equation*}
and their cyclic permutations. In each of these dimensions except for $(r,r,r,r)$, there is exactly one, up to isomorphism, indecomposable quadruple.

\section{Proof of Theorem \ref{tttt}}\label{ss-prof}

(a) Let $\mc T$ be a triple \eqref{as1a}. By Lemma \ref{jux}, $\mc T$ is isomorphic to some direct sum \eqref{ser03}. Write this sum in the form $\mc T_1\oplus \mc T_2$, in which $\mc T_1$ is the direct sum of all regular summands and $\mc T_2$ is the direct sum of the remaining summands. Let $\mc T_3$ be a strictly regular triple that is isomorphic to $\mc T_1$. Then the isomorphism $\mc T\simeq\mc T_3\oplus \mc T_2$ satisfies the statement (a) of Theorem \ref{tttt}.
\medskip

(b) This statement of Theorem \ref{tttt}
 is obvious.
\medskip

(c)
Let $\mc T$ be a strictly singular triple \eqref{as1a}. Using $\ind(\underline G)$ given in Lemma \ref{vbe2} and Steps 1--3 from Section \ref{red}, we can construct a direct sum \eqref{ser03}, which is  isomorphic to $\mc T$  by Lemma \ref{jux}.
Since $\mc T$ is strictly singular,
all summands of this direct sum are strictly singular. Thus, they cannot be obtained from the first quadruple in \eqref{las} with nonsingular $\Phi$ and from the last two quadruples  in \eqref{las}. Denote by $\ind'(\underline G)$ the set of all remaining quadruples \eqref{las1}--\eqref{las} and the quadruples dual to \eqref{las2}.

Let us apply  Steps 1--3 from Section \ref{red} to $\ind'(\underline G)$.

 \begin{description}
   \item[Step 1:] In this step,
   we construct a partition of
$\ind' (\underline{G})$
into 3 subsets as in \eqref{4.8d}:
\begin{equation}\label{8d}
{\ind' (\underline{G})}
=
\begin{tabular}{|c|c|}
  \hline
\multicolumn{2}{|c|}
{${\cal P}={\cal P}^{\circ}\vphantom{{\hat{N}}}$}\\
 \hline &\\[-12pt]
  $\; {\cal Q}\;  $&
  ${\cal Q}^{\circ}\text{ if
  ${\cal Q}^{\circ}\not
  \simeq{\cal Q}$}$
  \\  \hline
\end{tabular}\,,\ \
\begin{matrix}
 {\cal P}\in
\ind'_0(\underline{G}),\\[1pt]
{\cal
Q}\in\ind'_1(\underline{G}).
\end{matrix}
\end{equation}
 \end{description}

Let us prove that
\begin{itemize}
  \item
$\ind'_0(\underline G)$ can be taken
consisting of the quadruples
\begin{equation}\label{c1}
\begin{split}
\syma{r}{r}{r-1}{r-1}
{F_{r}}{F_{r}^{\T}
\ \begin{smallmatrix}
(\text{$r$ is even if $\varepsilon =-$,}\\
\text{\  $r$ is odd if $\delta =-$})
\end{smallmatrix}
}{Z_{r,\varepsilon}}
{Z_{r-1,\delta}}
\quad
\syma{r-1}{r-1}{r}{r}
{F_{r}^{\T}}{F_{r}
\ \begin{smallmatrix}
(\text{$r$ is odd if $\varepsilon =-$,}\\
\text{\quad $r$ is even if $\delta =-$})
\end{smallmatrix}
}{Z_{r-1,\varepsilon} }
{Z_{r,\delta} }
\end{split}
\end{equation}
in which $Z_{r,+}$ and $Z_{r,-} $ are defined in \eqref{asdq}, and

  \item
$\ind'_1(\underline G)$ can be taken
consisting of
\begin{itemize}

  \item[$(1^{\circ})$]
the first quadruple in \eqref{las1}, in which $r$ is odd and $\varepsilon =-$, or $r$ is even and $\delta=-$,

  \item[$(2^{\circ})$]
the second quadruple in \eqref{las1},
in which $r$ is even and $\varepsilon=-$, or $r$ is odd and $\delta=-$,

  \item[$(3^{\circ})$]
the quadruples \eqref{las2},

  \item[$(4^{\circ})$]
the first quadruple in \eqref{las} with $\Phi_r=J_r(0)^{\T}$.
\end{itemize}
\end{itemize}

The horizontal arrows in the quadruples \eqref{las1} are assigned by nonsingular matrices. Since every nonsingular skew-symmetric matrix is of even size, the  quadruples ($1^{\circ}$) and ($2^{\circ}$) are not isomorphic to self-dual quadruples.
The remaining quadruples \eqref{las1} are isomorphic to the self-dual quadruples \eqref{c1}. For example, if $r$ is odd and $\varepsilon =+$, then the first quadruple in \eqref{las1} is isomorphic to the first self-dual quadruple in \eqref{c1}
since $Z_{r}G^{\T}_rZ_{r-1}=F_r^{\T}$ and
\[
\syma{r}{r}{r-1}{r-1}{F_r}{G_r^{\T}}{I_r}
{I_{r-1}}
             \quad
{\xymatrix@R=11pt@C=50pt{\\
{}\ar@{|->}[r]
^{\left(\begin{smallmatrix}
               I_r& Z_r \\
               I_{r-1} & Z_{r-1}
             \end{smallmatrix}\right)}&{}}}
            \quad
\syma{r}{r}{r-1}{r-1}{F_r}{F_r^{\T}}
{Z_r}{Z_{r-1}}
\quad\xymatrix@=3pt{{}\\ \simeq}\quad
\symaa{r}{r}{r-1}{r-1}
{F_{r}}{F_{r}^{\T}}{Z_{r}}
{Z_{r-1,\delta}}
\]

The quadruples ($3^{\circ}$) are not isomorphic to self-dual quadruples since
if some quadruple \eqref{vrj} is isomorphic to a self-dual quadruple, then $m_1=m_2$ and $n_1=n_2$.

The quadruple ($4^{\circ}$) is dual to the second quadruple in \eqref{las}, and so it is not isomorphic to a self-dual quadruple.

\begin{description}
  \item[Steps 2:] Let us construct a set of triples
$\{\mc P^a\,|\,0\ne a=a^{\circ}\in \TT(\mc P)\}$
for each $\mc P\in\ind'_0(\underline G)$.
\end{description}
Let $P$ be the first quadruple in \eqref{c1} with odd $r$ and $\varepsilon =+$. Let $\mattc{R_1&S_1\\R_2&S_2}:P\to P$ be an endomorphism of $P$. Then
\[
F_rR_1=R_2F_r,\quad S_1F_r^{\T}=F_r^{\T}S_2,\quad
S_1Z_r=Z_rR_1,\quad
S_2Z_{r-1,\delta }=Z_{r-1,\delta } R_2.
\]
A straightforward  computation shows that
\[
(R_1,R_2,S_1,S_2)=\left(
\matt{a&&0\\[-5pt]&\ddots\\ *&&a}, \matt{a&&0\\[-5pt]&\ddots\\ *&&a}, \matt{a&&*\\[-5pt]&\ddots\\ 0&&a}, \matt{a&&*\\[-5pt]&\ddots\\ 0&&a}\right),\quad a\in\FF.
\]
Hence, we can identify $\TT(P)$ and $\FF$. The mapping $a\mapsto \phi_a:= \left(\begin{smallmatrix}aI_r&aI_r
\\aI_{r-1}&aI_{r-1}
\end{smallmatrix}\right)$ is an embedding of $\FF$ into $\End(P)$. If $0\ne a=\tilde a\in\FF$, then
\[
\symaaa
{r}{r}{r-1}{r-1}
{\smash{\textstyle P^{\phi_a}:}
\qquad F_{r}}{F_{r}^{\T}}{aZ_{r}}
{aZ_{r-1,\delta}}
\]
(see \eqref{ebh}) is a self-dual quadruple and the corresponding triple $P^a$ is the first triple in \eqref{b1a} with odd $r$ and $\varepsilon =+$.

In the same manner,
we can identify $\TT(P)$ and $\FF$ for each quadruple $P$ from \eqref{c1}. The mapping $a\mapsto \phi_a$, in which $ \phi_a:=\mattc{aI_r&aI_r
\\aI_{r-1}&aI_{r-1}}$ for the first quadruple in \eqref{c1} and $ \phi_a:=\mattc{aI_{r-1}&aI_{r-1}
\\aI_r&aI_r}$ for the second quadruple in  \eqref{c1},  is an embedding of $\FF$ into $\End(P)$. If $0\ne a=\tilde a\in\FF$, then the quadruple $P^{\phi_a}$ is obtained from $P$ by multiplying by $a$ the matrices that correspond to the horizontal arrows. The triple $P^{a}$ is the first triple in \eqref{b1a} or the first triple in \eqref{b2}.

\begin{description}
  \item[Step 3:] Let us construct the triple $\mc Q^+$ for each $\mc Q\in\ind_1'(\underline{G})$.
\end{description}
By \eqref{ekx} and the correspondence \eqref{ekx1}, if $\mc Q$ is ($1^{\circ}$) or ($2^{\circ}$), then $\mc Q^+$  is the second triple in \eqref{b1a} or \eqref{b2}, respectively.
The quadruples ($3^{\circ}$) give the triples \eqref{b5} and the first two triples in \eqref{b6}. The quadruple ($4^{\circ}$) gives the last triple in \eqref{b6}.

This proves the statement (c) due to Lemma \ref{jux}.

\medskip

(d) This statement follows from (c) since by the law of inertia
\begin{itemize}
\item each symmetric form over $\CC$ and each Hermitian form over $\HH$ with involution that differs from  \eqref{tgh} are reduced to exactly one form $\tilde x_1x_1+\dots+
\tilde x_kx_k$, and

  \item each symmetric form over $\RR$, each Hermitian form over $\CC$, and each Hermitian form over $\HH$ with involution \eqref{tgh} are reduced to exactly one form $\tilde x_1x_1+\dots+
\tilde x_lx_l-\tilde x_{l+1}x_{l+1}-\dots-\tilde x_kx_k$,
  \end{itemize}
in which the involution $a\mapsto\tilde a$ is the identity if the form is symmetric; see \cite[p.\,484]{ser_izv} for Hermitian forms over $\HH$. This proves Theorem \ref{tttt}.

\section*{Acknowledgements}
This paper is a result of a student seminar held at the University of S\~ao Paulo during the visit of V.V.~Sergeichuk in 2016;
he is grateful to the university for hospitality and to the FAPESP
for financial support (grant 2015/05864-9).

\end{document}